\documentclass[12pt]{article}
\usepackage{amsmath,amsfonts,amscd,a4wide,theorem}
\usepackage{pgf,pgfarrows, pgfnodes}
\theoremheaderfont{\bfseries\rmfamily\upshape}
                \newtheorem{theorem}{Theorem}
                \newtheorem{lemma}{Lemma}

{\theorembodyfont{\rmfamily}

                \newtheorem{remark}{Remark}

}

\newcommand{\trace}{\operatorname{trace}}

\newcommand{\hnabla}{\widehat{\nabla}}

\newcommand{\qed}{\text{\rule{.4em}{1.7ex}\hspace{.6em}}}
\newenvironment{proof}[1][]{\noindent {\bf Proof#1:\ }}
        {\hspace*{.1em}\hfill\qed\bigskip \noindent}
\newcounter{rom}
\renewcommand{\therom}{(\roman{rom})}
\newenvironment{romanlist}{\begin{list}{\therom}
                {\setlength{\leftmargin}{2em}\usecounter{rom}}}%
{\end{list}}
\title{A characterisation of the Calabi product of hyperbolic affine spheres}
\begin{document}
\author{Zejun Hu
\thanks {Supported by grants of NSFC-10671181 and Chinese-German
cooperation projects DFG PI 158/4-5.}
\and Haizhong Li\thanks{Supported by grants of NSFC-10531090 and
Chinese-German cooperation projects DFG PI 158/4-5.}
\and Luc Vrancken}
\date{}

\maketitle

\sloppy

\begin{abstract}\noindent There exists a well known construction which allows
to associate with two hyperbolic affine spheres $f_i: M_i^{n_i} \rightarrow
\mathbb R^{n_i+1}$ a new hyperbolic affine sphere immersion
of  $I
\times M_1 \times M_2$ into $\mathbb R^{n_1+n_2+3}$.  In
this paper we deal with the inverse problem: how to
determine from properties of the difference tensor whether a
given hyperbolic affine sphere immersion of a manifold $M^n \rightarrow \mathbb R^{n+1}$ can be
decomposed in such a way.
\end{abstract}

{{\bfseries Key words}:  {\em affine hypersphere, Calabi product, affine hypersurface}.

{\bfseries Subject class: } 53A15.}

\section{Introduction}

In this paper we study nondegenerate affine hypersurfaces $M^n$ into
$\mathbb R^{n+1}$, equipped with its standard affine connection $D$.
It is well known that on such a hypersurface there exists a
canonical transversal vector field $\xi$, which is called the affine
normal. With respect to this transversal vector field one can
decompose
\begin{equation}
D_X Y = \nabla_X Y +h(X,Y) \xi,
\end{equation}
thus introducing the affine metric $h$ and the induced affine
connection $\nabla$. The Pick-Berwald theorem states that $\nabla$
coincides with the Levi Civita connection $\hnabla$ of the affine
metric $h$ if and only if $M$ is immersed as a nondegenerate
quadric. The difference tensor $K$ is introduced by
\begin{equation}
K_X Y = \nabla_X Y -\hnabla_X Y.
\end{equation}
It follows easily that $h(K(X,Y),Z)$ is symmetric in $X$, $Y$ and
$Z$. The apolarity condition states that $\trace K_X =0$ for every
vector field $X$. The fundamental theorem of affine differential
geometry, Dillen, see Ref.~\cite{dinovr91} implies that an affine
hypersurface is completely determined by the metric and the
difference tensor $K$.

Deriving the affine normal, we introduce the affine shape operator $S$ by
\begin{equation}
D_X \xi =-SX.
\end{equation}

Here, we will restrict ourselves to the  case that the affine shape
operator $S$ is a multiple of the identity, i.e. $S=H I$. This means
that all affine normals are parallel or pass through a fixed point.
We will also assume that the metric is positive definite in which
case one distinguishes the following classes of affine hyperspheres:
\begin{romanlist}
\item elliptic affine hyperspheres, i.e. all affine normals pass through a fixed point and $H >0$,
\item hyperbolic affine hyperspheres, i.e. all affine normals pass through a fixed point and $H<0$,
\item parabolic affine hyperspheres, i.e. all the affine normals are parallel ($H=0$).
\end{romanlist}
Due to the work of amongst others Calabi \cite{ca72}, Pogorelov
\cite{po72}, Cheng and Yau \cite{chya86}, Sasaki \cite{sa80} and Li
\cite{li92}, positive definite affine hyperspheres which are
complete with respect to the affine metric $h$ are now well
understood. In particular, the only  complete elliptic or parabolic
positive definite affine hyperspheres are respectively the ellipsoid
and the paraboloid. However, there exist many hyperbolic affine
hyperspheres.

In the local case, one is far from obtaining a classification. The reason for this is that affine hyperspheres
reduce to the study of the Monge-Amp{\`e}re equations.
Calabi introduced a construction, called the Calabi product, which shows how to associate with one (or two)
hyperbolic affine hyperspheres a new hyperbolic affine hypersphere. This construction,
as well as the corresponding properties for the difference tensor are recalled in the next section.

In this paper we are interested in the reverse construction, i.e. how to determine using properties of the difference
tensor whether or not a given hyperbolic affine hypersphere (with mean curvature $-1$) can be decomposed
as a Calabi product of a hyperbolic affine hypersphere and a point or as a Calabi product of two hyperbolic
affine hyperspheres.

In particular we show the following two theorems:
\begin{theorem} Let $\phi: M^n \rightarrow \mathbb R^{n+1}$ be a (positive definite) hyperbolic affine hypersphere with mean
curvature $\lambda$, $\lambda<0$. Assume that there exists two
distributions $\mathcal D_1$ and $\mathcal D_2$ such that
\begin{romanlist}
\item $T_pM = \mathcal D_1 \oplus \mathcal D_2$,
\item $\mathcal D_1$ and $\mathcal D_2$ are orthogonal with respect to the affine metric $h$
\item $\mathcal D_1$ is a one dimensional distribution spanned by a unit length vector field $T$
\item there exist numbers $\lambda_1$ and $\lambda_2$ satisfying $-\lambda+\lambda_1 \lambda_2 -\lambda_2^2= 0$ such that
\begin{align*}
& K(T,T)=\lambda_1 T\\
&K(T,U) = \lambda_2 U,
\end{align*}
where $U \in \mathcal D_2$.
\end{romanlist}
Then $\phi:M^n \rightarrow \mathbb R^{n+1}$ can be decomposed as the Calabi product of a hyperbolic
affine sphere $\psi:M_1^{n-1} \rightarrow \mathbb R^{n}$ and a point.
\end{theorem}
and
\begin{theorem} \label{theoremprod2}Let $\phi: M^n \rightarrow \mathbb R^{n+1}$ be a (positive definite) hyperbolic affine hypersphere with mean
curvature $\lambda$, $\lambda<0$. Assume that there exists
distributions $\mathcal D_1$ (of dimension 1, spanned by a unit
length vector field $T$), $\mathcal D_2$
 (of dimension $n_2$) and $\mathcal D_3$ (of dimension $n_3$) such that
\begin{romanlist}
\item $1+n_2+n_3 = n$,
\item $\mathcal D_1$, $\mathcal D_2$  and $\mathcal D_3$ are mutually orthogonal
 with respect to the affine metric $h$
\item there exist numbers $\lambda_1$, $\lambda_2$ and $\lambda_3$ such that
\begin{align*}
& K(T,T)=\lambda_1 T\\
&K(T,V) = \lambda_2 V,\\
&K(T,W)= \lambda_3 W,\\
&K(V,W)=0.
\end{align*}
where $V \in \mathcal D_2$, $W \in \mathcal D_3$, $\lambda_1 =
\lambda_2 +\lambda_3$ and $\lambda_2 \lambda_3 = \lambda$.
\end{romanlist}
Then $\phi:M^n \rightarrow \mathbb R^{n+1}$ can be decomposed as the Calabi product of two hyperbolic
affine sphere immersions $\psi_1:M_1^{n_2} \rightarrow \mathbb R^{n_2+1}$ and
 $\psi_2:M_2^{n_3} \rightarrow \mathbb R^{n_3+1}$.
\end{theorem}
Note that, as explained in the next section, the converse of the above two theorems is also true.

To conclude this introduction, we remark that the basic integrability conditions for a hyperbolic
 affine hypersphere with mean curvature $-1$ state that:
\begin{align}
&\hat R(X,Y)Z = -(h(Y,Z) X-h(X,Z)Y)-[K_X,K_Y]Z,\label{gauss}\\
&(\hat\nabla K)(X, Y,Z) =(\hat \nabla K)(Y,X,Z).\label{codazzi}
\end{align}


\section{The Calabi product}

Let $\psi_1: M_1^{n_2} \rightarrow R^{n_2+1}$ and $\psi_2: M_2^{n_3} \rightarrow R^{n_3+1}$
be hyperbolic affine hyperspheres with mean curvature $-1$. Then we define
the Calabi product of $M_1$ with a point by
\begin{equation*}
\tilde \psi(t,p)= (c_1 e^{\tfrac{t}{\sqrt{n}}} \psi_1(p), c_2
e^{-{\sqrt{n}}t}),
\end{equation*}
where $p \in M_1$ and $t \in \mathbb R$
and the Calabi product of $M_1$ with $M_2$
by
\begin{equation*}
\psi(t,p,q))= (c_1 e^{\tfrac{\sqrt{n_3+1} t}{\sqrt{n_2+1}}}\psi_1(p),
 c_2 e^{-\tfrac{\sqrt{n_2+1} t}{\sqrt{n_3+1}}} \psi_2(q)),
\end{equation*}
where $p \in M_1$, $q \in M_2$ and $t \in \mathbb R$.

We now investigate the conditions on the constants $c_1$ and $c_2$ in order that the Calabi product
has constant mean curvature $-1$. We first do so for the Calabi product of two affine spheres.
We denote by $v_1,\dots,v_{n_2}$ local coordinates for $M_1$ and by $w_1,\dots,w_{n_3}$ local coordinates
for $M_2$. Then, it follows that
\begin{align*}
&\psi_t= (c_1 \tfrac{\sqrt{n_3+1}}{\sqrt{n_2+1}}e^{\tfrac{\sqrt{n_3+1} t}{\sqrt{n_2+1}}}\psi_1(p),
-c_2 \tfrac{\sqrt{n_2+1}}{\sqrt{n_3+1}}e^{-\tfrac{\sqrt{n_2+1} t}{\sqrt{n_3+1}}} \psi_2(q)),\\
&\psi_{tt}=(c_1 \tfrac{n_3+1}{n_2+1}e^{\tfrac{\sqrt{n_3+1} t}{\sqrt{n_2+1}}}\psi_1(p),
-c_2 \tfrac{n_2+1}{n_3+1}e^{-\tfrac{\sqrt{n_2+1} t}{\sqrt{n_3+1}}} \psi_2(q)),\\
&\psi_{tv_i}=\tfrac{\sqrt{n_3+1}}{\sqrt{n_2+1}}(c_1 e^{\tfrac{\sqrt{n_3+1} t}{\sqrt{n_2+1}}}(\psi_1)_{v_i},0),\\
&\psi_{tw_j}=-\tfrac{\sqrt{n_2+1}}{\sqrt{n_3+1}}(0,c_2 e^{-\tfrac{\sqrt{n_2+1} t}{\sqrt{n_3+1}}} (\psi_2)_{w_j}),\\
&\psi_{v_iv_j}=(c_1 e^{\tfrac{\sqrt{n_3+1} t}{\sqrt{n_2+1}}}(\psi_1)_{v_iv_j},0),\\
&\psi_{w_iw_j}=(0,c_2 e^{-\tfrac{\sqrt{n_2+1} t}{\sqrt{n_3+1}}} (\psi_2)_{w_iw_j}).
\end{align*}
If we denote by $h_2$ the affine metric on $M_2$ and by $h_3$ the
centroaffine metric introduced on $M_3$ it follows from the above
formulas that
\begin{align*}
&\psi_{tt} =\tfrac{n_3-n_2}{\sqrt{(n_2+1)(n_3+1)}} \psi_t + \psi\\
&\psi_{v_iv_j} = \tfrac{\sqrt{(n_2+1)(n_3+1)}}{n_2+n_3+2} h_2(\partial v_i,\partial v_j) \psi+...\\
&\psi_{w_iw_j} = \tfrac{\sqrt{(n_2+1)(n_3+1)}}{n_2+n_3+2}
h_3(\partial w_i,\partial w_j) \psi+...
\end{align*}
From \cite{nosa94} we see that $M$ is an affine hypersphere with
mean curvature $-1$ if and only if
\begin{equation*}
det[\psi,\psi_t,\psi_{v_1},\dots,\psi_{v_{n_2}},\psi_{w_1},\dots,\psi_{w_{n_3}}]^2
=h(\partial_t,\partial_t) det[h(\partial v_i,\partial v_j)]
det[h(\partial w_i,\partial w_j)].
\end{equation*}
Taking into account that $\psi_1$ and $\psi_2$ are already affine
spheres with mean curvature $-1$ we must have that
\begin{equation}
(c_1)^{n_2+1} (c_2)^{n_3+1} = \left
(\tfrac{\sqrt{(n_2+1)(n_3+1)}}{n_2+n_3+2}\right )^{n_2+n_3+2}.
\end{equation}
Hence we can take
\begin{align*}
&c_1=\tfrac{\sqrt{(n_2+1)(n_3+1)}}{n_2+n_3+2} d_1\\
&c_2=\tfrac{\sqrt{(n_2+1)(n_3+1)}}{n_2+n_3+2} d_2,
\end{align*}
where
\begin{equation}
(d_1)^{n_2+1} (d_2)^{n_3+1} = 1.
\end{equation}
Hence by applying an equiaffine transformation we may assume that
$d_1=d_2=1$ and therefore that the Calabi product of two hyperbolic
affine spheres with mean curvature $-1$ is an hyperbolic affine
sphere with mean curvature $-1$ if and only if
\begin{equation*}
\psi(t,p,q))= \tfrac{\sqrt{(n_2+1)(n_3+1)}}{n_2+n_3+2}(
e^{\tfrac{\sqrt{n_3+1} t}{\sqrt{n_2+1}}}\psi_1(p),
 e^{-\tfrac{\sqrt{n_2+1} t}{\sqrt{n_3+1}}} \psi_2(q)),
\end{equation*}
up to an equiaffine transformation.

For the Calabi product of a hyperbolic affine sphere and a point, we
proceed in the same way to deduce the following. The Calabi product
of a hyperbolic affine spheres with mean curvature $-1$ and a point
is an hyperbolic affine sphere with mean curvature $-1$ if and only
if
\begin{equation*}
\tilde \psi(t,p)= \tfrac{\sqrt{n}}{n+1}( e^{\tfrac{t}{\sqrt{n}}}
\psi_1(p), e^{-\sqrt{n} t}),
\end{equation*}
up to an equiaffine transformation.

\begin{remark} A straightforward calculation shows that the Calabi
product of two hyperbolic affine spheres has parallel cubic form
(with respect to the Levi Civita connection) if and only if both
original hyperbolic affine spheres have parallel cubic forms.
Similarly one has that the Calabi product of a hyperbolic affine
sphere and a point has parallel cubic form if and only if the
original affine sphere has parallel cubic form.
\end{remark}

\section{Characterisation of the Calabi product of two hyperbolic affine spheres and the proof of Theorem 2}

Throughout this section we will assume that $\phi:
M^n\longrightarrow\mathbb R^{n+1} $ is a hyperbolic affine
hypersphere. Without loss of generality we may assume that
$\lambda=-1$ by applying a homothety. We will now prove Theorem
\ref{theoremprod2}. Therefore, we shall also assume that $M$ admits
three mutually orthogonal differential distributions $\mathcal D_1$,
$\mathcal D_2$ and $\mathcal D_3$ of dimension $1$, $n_2> 0$ and
$n_3> 0$ respectively with $1+n_2+n_3=n$, and, for all vectors $V\in
\mathcal D_2$, $W\in \mathcal D_3$,
\begin{gather*}
K(T,T) = \lambda_1 T,\;\;\;\;\;\;
K(T,V) = \lambda_2 V,\\
K(T,W) = \lambda_3 W,\;\;\;\;\;\; K(V,W) = 0.
\end{gather*}
By the apolarity condition we must have that
\begin{equation}
\lambda_1 +n_2 \lambda_2 +n_3 \lambda_3 = 0,
\end{equation}
Moreover, we will assume that
\begin{align}
&\lambda_1 = \lambda_2+\lambda_3\\
&\lambda_2 \lambda_3 = -1.
\end{align}

The above conditions imply that $\lambda_1$, $\lambda_2$ and $\lambda_3$ are constants and
can be determined explicitly in terms of the dimension $n$.

As $M$ is a hyperbolic affine sphere we have that the difference
tensor is a symmetric tensor with respect to the Levi Civita
connection $\hat \nabla$ of the affine metric. In that case, as also
$h(K(X,Y),Z)$ is totally symmetric, the information of Lemma 1 and
Lemma 2 of \cite{BRV} remains valid and can be summarized in the
following lemma:
\begin{lemma} \label{lemme1}
We have
\begin{enumerate}
\item $\hat \nabla_{\mathcal D_1} \mathcal D_1 \subset \mathcal D_1$
\item $\hat \nabla_{\mathcal D_2} \mathcal D_2 \subset \mathcal D_2
\oplus \mathcal D_3$
\item $\hat \nabla_{\mathcal D_3} \mathcal D_3
\subset \mathcal D_2 \oplus \mathcal D_3$
\item $h(\hat \nabla_T
W,V) = h(\hat \nabla_W T,V)=-h(\hat \nabla_V T,W)$, for any $V
\in\mathcal D_2, W\in\mathcal D_3$
\end{enumerate}
\end{lemma}

Similarly using the
information of the previous lemma, Lemma 3 of \cite{BRV} reduces to
\begin{lemma} \label{lemme2} We have
\begin{enumerate}
\item $(\lambda_3-\lambda_2) h(\hat \nabla_{V} \tilde V, W)
= h( K(V,\tilde V),  \hat\nabla_T W)
 =  h(K(V,\tilde V),  \hat \nabla_W T)$
\item  $(\lambda_2-\lambda_3) h(\hat\nabla_{W} \tilde W, V)
=  h(K(W,\tilde W),  \hat\nabla_T V)
 =  h(K(W,\tilde W), \hat \nabla_V T)$
\end{enumerate}
\end{lemma}

We denote now by $\{V_1,\dots,V_{n_2}\}$, respectively  $\{W_1,\dots,W_{n_3}\}$ an
orthonormal basis of $\mathcal D_2$ (resp. $\mathcal D_3$) with respect to the affine
metric $h$. Then, we have
\begin{lemma} \label{lemme3} Let $V, \tilde V \in \mathcal D_2$. Then
$$h(\hat \nabla_V T, \hat \nabla_{\tilde V} T) =0.$$
\end{lemma}
\begin{proof} Using the Gauss equation, we have that
\begin{align*}
h(\hat R(V,T)T,\tilde V)&= -h(V,\tilde V) - h(K(T,T),K(V,\tilde V))
 + h(K(T,V),K(T,\tilde V))\\
&= (-1- \lambda_1 \lambda_2 +\lambda_2^2) h( V,\tilde V)\\
&=(-1- \lambda_3 \lambda_2) h(V,\tilde V)=0.
\end{align*}
On the other hand, by a direct computation using the previous lemmas, we have
\begin{align*}
h(\hat R(V,T)T,\tilde V)&=h(\hat \nabla_{V} \hat \nabla_{T}T -\hat
\nabla_{T}\hat \nabla_{V} T
 -\hat \nabla_{\hat \nabla_{V}T -\hat \nabla_{T}V} T,\tilde V)\\
&=h(-\hat\nabla_{T} \hat\nabla_V T,\tilde V ) -\sum_{k=1}^{n_3} h(\hat \nabla_{V}T - \hat \nabla_T V,W_k) h(\hat
\nabla_{W_k} T, \tilde V)\\
&=h(-\hat\nabla_{T} \hat \nabla_V T,\tilde V) \qquad \text{by Lemma \ref{lemme1} (iv)}\\
&=-\sum_{k=1}^{n_3} h(\hat\nabla_V T,W_k) h(\hat \nabla_{T} W_k,\tilde V)\\
&=\sum_{k=1}^{n_3} h(\hat \nabla_V T,W_k) h(\hat \nabla_{\tilde V} T,W_k)\\
&=h(\hat \nabla_V T, \hat \nabla_{\tilde V} T).
\end{align*}
\end{proof}

Similarly, we have
\begin{lemma} \label{lemme4} Let $W, \tilde W \in \mathcal D_3$. Then
$$h(\hat\nabla_W T, \nabla_{\tilde W} T) = 0.$$
\end{lemma}

Combining the two previous lemmmas with  Lemma \ref{lemme2} and
Lemma \ref{lemme1} we see that the distributions determined by
$\mathcal D_2$ and $\mathcal D_3$ are totally geodesic. It also
implies that $h(\hat \nabla_V T,W) = h(\hat \nabla_W T, V) = 0$.

This is sufficient to conclude that locally $(M,h)$ is isometric with
$I \times M_1 \times M_2$ where $T$ is tangent to $I$, $\mathcal D_2$ is tangent
to $M_1$ and $\mathcal D_3$ is tangent to $M_2$.

The product structure of $M$ implies the existence of local
coordinates $(t,p,q)$ for $M$ based on an open subset containing
the origin of $\mathbb R \times \mathbb R^{n_2} \times \mathbb R^{n_3}$,
 such that $\mathcal D_1$ is given by $dp=dq=0$,
$\mathcal D_2$ is given by $dt=dq=0$, and $\mathcal D_3$ is given by
$dt=dp=0$. We may also assume that $T = \tfrac{\partial}{\partial
t}$. We now put

\begin{equation}\label{eqphi2phi3}
\quad \phi_2= -f \lambda_3 \phi + f T, \quad \phi_3 =
 g \lambda_2 \phi - g T,
\end{equation}
where the functions $f$ and  $g$, which depend only on the variable $t$,  are determined by
\begin{align*}
f'=f (\lambda_3  - \lambda_1) ,\\
g'=g(\lambda_2 - \lambda_1).
\end{align*}
It is clear that solutions are given by
\begin{equation*}
f(t) = d_1 e^{(\lambda_3-\lambda_1)t} \qquad \text{and} \qquad g(t)
= d_2 e^{(\lambda_2-\lambda_1)t},
\end{equation*}
where $d_1$ and $d_2$ are constants. Of course, as $\lambda_1 =
\lambda_2+\lambda_3$ we can rewrite the above equation as
\begin{equation*}
f(t) = d_1 e^{-\lambda_2 t} \qquad \text{and} \qquad g(t) = d_2
e^{-\lambda_3t}.
\end{equation*}
Computing $\lambda_1$, $\lambda_2$ and $\lambda_3$ explicitly, where if necessary by changing the sign of $E_1$ we
may assume that $\lambda_2 \ge 0$ we find that
\begin{align*}
&\lambda_2 = \tfrac{\sqrt{n_3+1}}{\sqrt{n_2+1}},\\
&\lambda_3 = -\tfrac{\sqrt{n_2+1}}{\sqrt{n_3+1}}.
\end{align*}
Solving now the above equations for the immersion $\phi$ we find that
\begin{align*}
\phi &=\tfrac{1}{f (\lambda_2-\lambda_3)}\phi_2 - \tfrac{1}{g (\lambda_2-\lambda_3)} \phi_3\\
&=(\tfrac{1}{d_1} e^{\tfrac{\sqrt{n_3+1}}{\sqrt{n_2+1}} t} \phi_2
 + \tfrac{1}{d_2} e^{-\tfrac{\sqrt{n_2+1}}{\sqrt{n_3+1}} t} \phi_3)(\frac{\sqrt{(n_2+1)(n_3+1)}}{n_2+n_3+2}).
\end{align*}

A straightforward computation, using \eqref{eqphi2phi3}, now shows that
\begin{align*}
D_T (\phi_2)&=D_T(-f \lambda_3 \phi + f T)\\
&=f(\lambda_3-\lambda_1) (-\lambda_3\phi+T)+ f (-\lambda_3 T +(K(T,T)+\phi))\\
&=f(\lambda_3-\lambda_1)(-\lambda_3\phi+T)+ f ((\lambda_1-\lambda_3) T +\phi))\\
&=f (\lambda_2 \lambda_3 +1)\phi =0.
\end{align*}
Similarly
\begin{align*}
&D_W (\phi_2)=f(- \lambda_3 W +K(W,T))=0,\\
&D_T (\phi_3)=0,\\
&D_V (\phi_3)=0.
\end{align*}
The above implies that $\phi_2$ reduces to a map of $M_1$ in $\mathbb R^n$ whereas
$\phi_3$ reduces to a map of $M_2$ in $\mathbb R^n$. As we have
that
\begin{align*}
&d\phi_2(V)=D_V (\phi_2)=f(- \lambda_3 V +K(V,T))=f(-\lambda_3+\lambda_2)V,\\
&d\phi_3(W)=D_W(\phi_3)=g(\lambda_2 W -K(W,T))=g
(\lambda_2-\lambda_3)W,
\end{align*}
these maps are actually immersions. Moreover, denoting by $\nabla^1$ the
$\mathcal D_2$ component of $\nabla$, we find that
\begin{align*}
D_{V}d\phi_2(\tilde V)&=f(-\lambda_3+\lambda_2)D_{V}\tilde V\\
&= f(-\lambda_3+\lambda_2)\nabla_{V}\tilde V+ f(-\lambda_3+\lambda_2)h(V,\tilde V)\phi \\
&= f(-\lambda_3+\lambda_2)\nabla^1_{V}\tilde V+f(-\lambda_3+\lambda_2)(h(K(V,\tilde V),T)T+h(V,\tilde V)\phi)\\
&=d\phi_2(\nabla^1_{V}  \tilde V)+ f(-\lambda_3+\lambda_2)h(V,\tilde V)(\lambda_2 T +\phi)\\
&=d\phi_2(\nabla^1_{V} \tilde V)+f (-\lambda_3+\lambda_2)\lambda_2 h(V,\tilde V)(T-\lambda_3 \phi)\\
&=d\phi_2(\nabla^1_{V} \tilde V)+ (-\lambda_3+\lambda_2)\lambda_2
h(V,\tilde V)\phi_2.
\end{align*}
The above formulas imply that $\phi_2$ can be interpreted as a
centroaffine immersion contained in an $n_2+1$-dimensional vector
subspace of $\mathbb R^{n+1}$ with induced connection $\nabla^1$ and
affine metric $h_1 = (-\lambda_3+\lambda_2)\lambda_2 h$. Similarly,
we get that $\phi_3$ can be interpreted as a centroaffine immersion
contained in an $n_3+1$-dimensional vector subspace of $\mathbb
R^{n+1}$ with induced connection $\nabla^2$ (the restriction of
$\nabla$ to $\mathcal D_3$) and affine metric $h_2 =g
(\lambda_3-\lambda_2)\lambda_3 h$. Of course as both spaces are
complementary, we may assume by a linear transformation that the
$n_2+1$ dimensional space is spanned by the first $n_2+1$
coordinates of $\mathbb R^{n+1}$ whereas the $n_3+1$ dimensional
space is spanned by the last $n_3+1$ coordinates of $\mathbb
R^{n+1}$.

Moreover, taking $V_1,\dots, V_{n_2}$ as before, we find that
\begin{align*}
\sum_{i=1}^{n_2} (\nabla^1 h_1) (V,V_i,V_i)&=
\lambda_2(\lambda_2-\lambda_3)\sum_{i=1}^{n_2} (\nabla^1 h)
(V,V_i,V_i)\\
&=-2\lambda_2(\lambda_2-\lambda_3)\sum_{i=1}^{n_2} h(\nabla^1_V V_i,V_i)\\
&=-2\lambda_2(\lambda_2-\lambda_3)\sum_{i=1}^{n_2} h(\nabla_V V_i,V_i)\\
&=\lambda_2(\lambda_2-\lambda_3)\sum_{i=1}^{n_2} (\nabla h) (V,V_i,V_i)=0,
\end{align*}
as by assumption $h(K(V,W),W))= h(K(V,T),T)$. So $M_1$ is an hyperbolic
affine hypersphere. Choosing now the constant $d_1$ appropriately we may
assume that $M_1$ has mean curvature $-1$. A similar argument also holds for $M_2$.

As
\begin{align*}
\phi &=\tfrac{1}{f (\lambda_2-\lambda_3)}\phi_2 - \tfrac{1}{g (\lambda_2-\lambda_3)} \phi_3\\
&=(\tfrac{1}{d_1} e^{\tfrac{\sqrt{n_3+1}}{\sqrt{n_2+1}} t} \phi_2
 + \tfrac{1}{d_2} e^{-\tfrac{\sqrt{n_2+1}}{\sqrt{n_3+1}} t} \phi_3)(\frac{\sqrt{(n_2+1)(n_3+1)}}{n_2+n_3+2}).
\end{align*}
We note from Section 2 that we must have that $d_1^{n_2+1} d_2^{n_2+1} = 1$ and that therefore $\phi$ is given
as the Calabi product of the immersions $\phi_1$ and $\phi_2$.

\begin{remark} In case that $M$ has parallel difference tensor, i.e. if $\hat \nabla K=0$, the conditions
of Theorem 2 can be weakened. Indeed we can prove:
\begin{theorem} Let $M$ be a hyperbolic affine sphere with mean curvature $\lambda$,
where $\lambda<0$. Suppose that $\hat \nabla  K= 0$
and there exists $h$-orthonormal distributions $\mathcal D_1$ (of dimension $1$),
$\mathcal D_2$ (of dimension $n_2$)
and such that  $\mathcal D_3$ (of dimension $n_3$) such that
\begin{align*}
&K(T,T)=\lambda_1 T,\\
&K(T,V)=\lambda_2 V,\\
&K(T,W)=\lambda_3 W,
\end{align*}
where $T$ is a unit vector spanning $\mathcal D_1$ and $V \in
\mathcal D_2$, $W \in \mathcal D_3$. Moreover we suppose that
$\lambda_2 \ne \lambda_3$ and $2\lambda_2 \ne \lambda_1 \ne 2
\lambda_3$. Then $\phi:M^n \rightarrow \mathbb R^{n+1}$ can be
decomposed as the Calabi product of two hyperbolic affine sphere
immersions $\psi_1:M_1^{n_2} \rightarrow \mathbb R^{n_2+1}$ and
 $\psi_2:M_2^{n_3} \rightarrow \mathbb R^{n_3+1}$ with parallel cubic form.
\end{theorem}
\begin{proof} By applying an homothety we may choose $\lambda=-1$.
As $\hat \nabla K=0$, we also have that $\hat R. K = 0$. This means that
\begin{equation*}
\hat R(X,Y)K(Z,U)= K(\hat R(X,Y)Z,U)+K(Z, \hat R(X,Y)U).
\end{equation*}
So, taking $X=Z=U=T$ and $Y=V$, we find that
$$
\hat R(T,V)T=V -K_T K_V T +K_V K_TT=(1-\lambda_2^2+\lambda_1 \lambda_2) V.
$$
Hence we deduce that
\begin{equation*}
(\lambda_1-2\lambda_2)(-1-\lambda_1 \lambda_2 +\lambda_2^2)=0.
\end{equation*}
Similarly we have
\begin{equation*}
(\lambda_1-2\lambda_3)(-1-\lambda_1 \lambda_3 +\lambda_3^2)=0.
\end{equation*}
In view of the conditions, we must have that $\lambda_2$ and $\lambda_3$ are the two different
roots of the equation
$$-1 - \lambda_1 x +x^2=0.$$
Consequently $\lambda_2+\lambda_3 = \lambda_1$ and $\lambda_2 \lambda_3=-1$.

Finally we take $Z=U=T$, $X=V$ and $Y=W$. Then we find that
\begin{align*}
\lambda_1 \hat R(V,W)T&= 2 K(\hat R(V,W)T,T)=-2  K(K_V K_W T,T)+2
K(K_W K_V T,T)\\
&=-2(\lambda_3-\lambda_2)K_T K_VW.
\end{align*}
Hence
\begin{align*}
\lambda_1 (\lambda_2-\lambda_3) K_V W&= 2 K(\hat R(V,W)T,T)=-2 K(K_V
K_W T,T)+2 K(K_W K_V T,T)\\
&=-2(\lambda_3-\lambda_2)K_T K_VW.
\end{align*}
This implies that $K_V W$ is an eigenvector of $K_T$ with eigenvalue $\tfrac{1}{2} \lambda_1$. Given the form
of $K_T$ we deduce that $K(V,W) = 0$. We are now in a position to apply Theorem 2 and deduce that
$M$ can be obtained as the Calabi product of the hyperbolic affine spheres.
\end{proof}

\end{remark}

\section{Characterisation of the Calabi product of a hyperbolic affine sphere and a point and the proof of Theorem 1}

Throughout this section we will assume that $\phi:
M^n\longrightarrow\mathbb R^{n+1} $ is a hyperbolic affine
hypersphere with mean curvature $-1$ and we will prove Theorem 1.
Therefore, we shall also assume that $M$ admits two mutually
orthogonal differential distributions $\mathcal D_1$ and $\mathcal
D_2$ of dimension $1$ and $n_2> 0$, respectively, with $1+n_2=n$,
and, for unit vector $T\in \mathcal D_1$ and all vectors $V\in
\mathcal D_2$,
\begin{gather*}
K(T,T) = \lambda_1 T,\;\;\;\;\;\;
K(T,V) = \lambda_2 V.\\
\end{gather*}
By the apolarity condition we must have that
\begin{equation}
\lambda_1 +n_2 \lambda_2= 0,
\end{equation}
Moreover, we will assume that
\begin{equation}
1+\lambda_1 \lambda_2-\lambda_2^2=0.
\end{equation}

The above conditions imply that $\lambda_1$ and $\lambda_2$ are constant and
can be determined explicitly in terms of the dimension $n$.
Indeed, if necessary by replacing $T$ with $-T$, we have that
\begin{align*}
&\lambda_2= \tfrac{1}{\sqrt{n}},\\
&\lambda_1=-\tfrac{n-1}{\sqrt{n}}.
\end{align*}

We now proceed as in the previous case. Using the fact that
$\hat \nabla K$ is totally symmetric it follows that
\begin{lemma} \label{lemme3} We have
\begin{enumerate}
\item $\hat \nabla_{T} T=0$,
\item $\hat \nabla_{V} T=0$,
\item $h(\hat \nabla_V \tilde V,T)=0$.
\end{enumerate}
\end{lemma}

The  previous lemmma tells us that the distributions determined by
$\mathcal D_1$ and $\mathcal D_2$ are totally geodesic.
This is sufficient to conclude that locally $(M,h)$ is isometric with
$I \times M_1$ where $T$ is tangent to $I$ and $\mathcal D_2$ is tangent
to $M_1$.

The product structure of $M$ implies the existence of local
coordinates $(t,p$ for $M$ based on an open subset containing
the origin of $\mathbb R \times \mathbb R^{n_2}$,
 such that $\mathcal D_1$ is given by $dp=0$ and
$\mathcal D_2$ is given by $dt=$. We may also assume that $T = \tfrac{\partial}{\partial
t}$. We now put

\begin{equation}\label{eqphi2}
\quad \phi_2= f  \tfrac{1}{\lambda_2} \phi + f T, \quad \phi_3 =
 g \lambda_2 \phi - g T,
\end{equation}
where the functions $f$ and  $g$, which depend only on the variable $t$,  are determined by
\begin{align*}
f'=-f \lambda_2= -\tfrac{1}{\sqrt{n}} ,\\
g'=g(\lambda_2 - \lambda_1)=\sqrt{n}.
\end{align*}
It is clear that solutions are given by
\begin{equation*}
f(t) = d_1 e^{-\tfrac{1}{\sqrt{n}}t} \qquad \text{and} \qquad g(t)
= d_2 e^{\sqrt{n} t}.
\end{equation*}

A straightforward computation, now shows that
\begin{align*}
D_T (\phi_2)&=D_T(f \sqrt{n} \phi + f T)\\
&=-f(\phi+ \tfrac{1}{\sqrt{n}}T)+ f ( \sqrt{n} T +(K(T,T)+\phi))\\
&=f T (- \tfrac{1}{\sqrt{n}}+ \sqrt{n}-\tfrac{n-1}{\sqrt{n}})\\
&=0.
\end{align*}
Similarly
\begin{align*}
&D_T (\phi_3)=0,\\
&D_V (\phi_3)=0.
\end{align*}
The above implies that $\phi_2$ reduces to a map of $M_1$ in $\mathbb R^n$ whereas
$\phi_3$ is a constant vector in $\mathbb R^n$. As we have
that
\begin{equation*}
d\phi_2(V)=D_V (\phi_2)=f(\sqrt{n} V
+K(V,T))=f(\sqrt{n}+\tfrac{1}{\sqrt{n}})V,
\end{equation*}
the map $\phi_2$ is actually immersions. Moreover, denoting by $\nabla^1$ the
$\mathcal D_2$ component of $\nabla$, we find that
\begin{align*}
D_{V}d\phi_2(\tilde V)&=f(\sqrt{n}+\tfrac{1}{\sqrt{n}})D_{V}\tilde V\\
&= f(\sqrt{n}+\tfrac{1}{\sqrt{n}})\nabla_{V}\tilde V+ f(\sqrt{n}+\tfrac{1}{\sqrt{n}})h(V,\tilde V)\phi \\
&= f(\sqrt{n}+\tfrac{1}{\sqrt{n}})\nabla^1_{V}\tilde V+f(\sqrt{n}+\tfrac{1}{\sqrt{n}})(h(K(V,\tilde V),T)T+h(V,\tilde V)\phi)\\
&=d\phi_2(\nabla^1_{V}  \tilde V)+ f(\sqrt{n}+\tfrac{1}{\sqrt{n}})h(V,\tilde V)(\tfrac{1}{\sqrt{n}} T +\phi)\\
&=d\phi_2(\nabla^1_{V} \tilde V)+\tfrac{n+1}{n}
h(V,\tilde V)\phi_2.
\end{align*}
The above formulas imply that $\phi_2$ can be interpreted as a
centroaffine immersion contained in an $n_2+1$-dimensional vector
subspace of $\mathbb R^{n+1}$ with induced connection $\nabla^1$ and
affine metric $h_1 = \tfrac{n+1}{n} h$. Of course as the vector $\phi_3$ is transversal
to the immersion $\phi_2$, we may assume by a linear transformation that the
$\phi_2$ lies in the space spanned by the first $n$
coordinates of $\mathbb R^{n+1}$ whereas the constant vector lies in the direction of the last coordinate,
and by choosing $d_2$ appropriately we may assume that $\phi_2= (0,\dots,0,1)$.

As before we get that $M_1$ satisfies the apolarity condition and hence is a hyperbolic
affine hypersphere. Choosing now the constant $d_1$ appropriately we may
assume that $M_1$ has mean curvature $-1$.

As
\begin{equation*}
\phi
=(\tfrac{1}{d_1} e^{\tfrac{1}{\sqrt{n}} t} \phi_2
 + \tfrac{1}{d_2} e^{-\sqrt{n} t} \phi_3)(\frac{\sqrt{(n)(n_3+1)}}{n+1}).
\end{equation*}
We note from Section 2 that we must have that $d_1^{n_2+1} d_2 = 1$ and that therefore $\phi$ is given
as the Calabi product of the immersions $\phi_1$ and a point.


\vskip 1cm
\begin{flushleft}
\medskip\noindent
Zejun Hu: {\sc Department of Mathematics, Zhengzhou University,
Zhengzhou 450052, People's Republic of China.}\ \ E-mail:
huzj@zzu.edu.cn

\medskip
Haizhong Li: {\sc Department of Mathematical Sciences, Tsinghua
University, Beijing 100084, People's Republic of China} \ \ E-mail:
hli@math.tsinghua.edu.cn

\medskip
Luc Vrancken: {\sc LAMATH, ISTV2, Campus du mont houy, Universite de
Valenciennes, France.} \ \ E-mail: luc.vrancken@univ-valenciennes.fr

\end{flushleft}

\end{document}